\documentclass[11pt,a4paper]{amsart}

\usepackage{graphicx,amsthm,amsmath,amsfonts,calc,enumerate}
\usepackage[numbers,sort&compress]{natbib}
\usepackage[lmargin=30mm,rmargin=30mm,tmargin=30mm,bmargin=30mm,footskip=5mm]{geometry}

\renewcommand{\baselinestretch}{1.2}
\newcommand{\Figure}[4][htb]{
\begin{figure}[#1]
	\vspace*{1ex}
	\begin{center}#3\end{center}
	\vspace*{-1ex}
	\caption{\figlabel{#2}#4}
\end{figure}
}

\theoremstyle{plain}

\newtheorem{lemma}{Lemma}

\theoremstyle{definition}

\newcommand{\seclabel}[1]{\label{sec:#1}}
\newcommand{\secref}[1]{Section~\ref{sec:#1}}

\newcommand{\lemlabel}[1]{\label{lem:#1}}
\newcommand{\lemref}[1]{Lemma~\ref{lem:#1}}
\newcommand{\twolemref}[2]{Lemmas~\ref{lem:#1} and \ref{lem:#2}}

\newcommand{\figlabel}[1]{\label{fig:#1}}
\newcommand{\figref}[1]{Figure~\ref{fig:#1}}

\newcommand{\eqnlabel}[1]{\label{eqn:#1}}
\newcommand{\eqnref}[1]{\eqref{eqn:#1}}
\newcommand{\Eqnref}[1]{Equation~\eqref{eqn:#1}}

\newcommand{\half}{\ensuremath{\protect\tfrac{1}{2}}}
\newcommand{\bracket}[1]{\ensuremath{\protect\left(#1\right)}}
\newcommand{\ceil}[1]{\ensuremath{\protect\left\lceil #1\right\rceil}}
\newcommand{\arc}[1]{\ensuremath{\protect\overrightarrow{#1}}}
\newcommand{\Oh}[1]{\ensuremath{\protect\mathcal{O}(#1)}}
\newcommand{\cn}[1]{\ensuremath{\chi(#1)}}
\newcommand{\hcn}[1]{\ensuremath{\textsf{h}(#1)}}

\newcommand{\ocn}[1]{\ensuremath{\overrightarrow{\chi}(#1)}}

\begin{document}

\title[On the Oriented Chromatic Number of Dense Graphs]{On the Oriented Chromatic Number\\ of Dense Graphs}
\author{David R.\ Wood}
\address{Departament de Matem{\`a}tica Aplicada II, Universitat Polit{\`e}cnica de Catalunya, Barcelona, Spain}
\email{david.wood@upc.edu}
\thanks{Supported by a Marie Curie Fellowship of the European Community under contract 023865, and by the projects MCYT-FEDER BFM2003-00368 and Gen.\ Cat 2001SGR00224.}

\begin{abstract}
Let $G$ be a graph with $n$ vertices, $m$ edges, average degree $\delta$, and maximum degree $\Delta$. The \emph{oriented chromatic number} of $G$ is the maximum, taken over all orientations of $G$, of the minimum number of colours in a proper vertex colouring such that between every pair of colour classes all edges have the same orientation. We investigate the oriented chromatic number of graphs, such as the hypercube, for which $\delta\geq\log n$. We prove that every such graph has oriented chromatic number at least $\Omega(\sqrt{n})$. In the case that $\delta\geq(2+\epsilon)\log n$, this lower bound is improved to $\Omega(\sqrt{m})$. Through a simple connection with harmonious colourings, we prove a general upper bound of $\Oh{\Delta\sqrt{n}}$ on the oriented chromatic number. Moreover this bound is best possible for certain graphs. These lower and upper bounds are particularly close when $G$ is ($c\log n$)-regular for some constant $c>2$, in which case the oriented chromatic number is between $\Omega(\sqrt{n\log n})$ and \Oh{\sqrt{n}\log n}.
\end{abstract}


\date{\today}

\maketitle

\section{Introduction}

Throughout this paper, $G$ is a (finite and simple) undirected graph with $n$ vertices, $m$ edges, and maximum degree $\Delta$. A \emph{colouring} of $G$ is a function that assigns a `colour' to each vertex so that adjacent vertices receive distinct colours. The \emph{chromatic number} \cn{G} is the minimum number of colours in a colouring of $G$. An \emph{oriented graph} $D$ is a directed graph with no parallel and no antiparallel arcs; that is, no two arcs have the same pair of endpoints. If $G$ is the underlying undirected graph of $D$ then $D$ is an \emph{orientation} of $G$. An \emph{oriented colouring} of $D$ is a colouring of $G$ such that between each pair of colour classes, all edges have the same direction. That is, there are no arcs \arc{vw} and \arc{xy} with $c(v)=c(y)$ and $c(w)=c(x)$. The \emph{oriented chromatic number} \ocn{D} is the minimum number of colours in an oriented colouring of $D$. The \emph{oriented chromatic number} \ocn{G} is the maximum of \ocn{D}, taken over all orientations $D$ of $G$. The oriented chromatic number was introduced by 
\citet{Courcelle-DAM94} in 1994 and is now a widely studied parameter; see \citep{Courcelle-DAM94, NesRasSop-DM97, KSZ-JGT97, BKNRS-DM98, BKNRS-DM99, BFKRS-JCTB01, Sopena-JGT97, EO-IPL06, Sopena-DM01, KlosMac-DM04, KlosMac-BICA04, RasSop-IPL94, Samal03, FHPZ-GC98, Ochem-IPL04, Wood-DMTCS05, DolSop-DM06, NesRas99, BorIva05, BI05, FRR-IPL03, ST-IPL04, Samal04, Sopena-IPL02, SS99, NSV97, KlosMac-DMGT04, KLSS99, HosSop-IPL06}.

This paper is motivated by a question of Andr\'e Raspaud [private communication, Prague 2004], who asked for the oriented chromatic number of the $d$-dimensional \emph{hypercube} $Q_d$. This is the graph with vertex set $\{0,1\}^d$, where two vertices are adjacent whenever they differ in precisely one coordinate. $Q_d$ is $d$-regular and has $2^d$ vertices. In this paper we prove generally applicable bounds on \ocn{G}, which in the case of the hypercube give 
\begin{equation}
\eqnlabel{Hypercube}
0.8007...\,\sqrt{2^d}\;\leq\;\ocn{Q_d}\;\leq\;2d\sqrt{2^d-1}\enspace,
\end{equation}
thus determining \ocn{Q_d} to within a factor of about $\frac{5}{2}d$. No non-trivial bounds on \ocn{Q_d} were previously known, as we now describe.

An undirected graph has $\chi(G)=n$ if and only if $G=K_n$; that is, if $G$ has diameter $1$. But when does $\ocn{G}=n$? This question was asked by \citet{FHPZ-GC98}, who observed that for every oriented graph $D$, 
\begin{equation}
\eqnlabel{DiameterTwo}
\text{$\ocn{G}=n$ if and only if $D$ has diameter $2$.}
\end{equation}
Here the \emph{diameter} of $D$ is the least integer $k$ such that every pair of vertices in $D$ are connected by a directed path of at most $k$ edges. 
\citet{KlosMac-DMGT04} called an oriented graph with diameter $2$ an \emph{oclique}. Note that small diameter ($>2$) does not necessarily imply large oriented chromatic number. For example, $K_{1,1,n}$ has an orientation with diameter $3$ and oriented chromatic number $3$. \citet{ERS66} proposed studying the extremal function $f(n)$, defined to be the minimum number of arcs in an oriented graph with $n$ vertices and diameter $2$. \citet{KS67} proved that $\frac{n}{2}\log\frac{n}{2}\leq f(n)\leq n\ceil{\log n}$, and \citet{FHPZ-GC98} tightened both bounds to conclude that $f(n)=(1-o(1))n\log n$. These results imply that there are $n$-vertex graphs with the same number of edges as the hypercube (that is, $n\log n$), yet have oriented chromatic number $n$. Thus good bounds for \ocn{Q_d} cannot be obtained just in terms of the number of edges. 

The example of an oriented graph with diameter $2$ by \citet{FHPZ-GC98} has a vertex of degree $n-1$. Thus it is natural to consider the oriented chromatic number of graphs with bounded degree.  \citet{Sopena-JGT97} and \citet{KSZ-JGT97} proved that the oriented chromatic number is bounded for graphs of bounded degree. The best bound is due to \citet{KSZ-JGT97}, who proved that every graph $G$ satisfies
\begin{equation}
\eqnlabel{DegreeUpperBound}
\ocn{G}\leq2\Delta^22^\Delta,
\end{equation}
and if $G$ is $\Delta$-regular with sufficiently many vertices then \begin{equation}
\eqnlabel{DegreeLowerBound}
\ocn{G}\geq2^{\Delta/2}.
\end{equation}
Thus the exponential dependence on $\Delta$ in \eqnref{DegreeUpperBound} is unavoidable. Observe that for graphs such as the hypercube with $\Delta\geq\log n$, the upper bound in \eqnref{DegreeUpperBound} is greater than the trivial upper bound of $n$. 

This motivates the study of the oriented chromatic number of graphs whose average degree is at least logarithmic in the number of vertices. For any such graph we establish a lower bound of $\Omega(\sqrt{n})$ on the oriented chromatic number. If the average degree is at least $(2+\epsilon)\log n$ then this lower bound is improved to $\Omega(\sqrt{m})$. These results are proved in \secref{UniversalLowerBound}. In \secref{UpperBound} we use a simple connection with harmonious colourings to prove a general upper bound of $\Oh{\Delta\sqrt{n}}$ on the oriented chromatic number. Moreover this bound is best possible for certain graphs, as proved in \secref{ExistentialLowerBound}. 

\section{An Upper Bound}
\seclabel{UpperBound}

In this section we prove an elementary upper bound on the oriented chromatic number. A colouring of an undirected graph $G$ is \emph{harmonious} if the endpoints of every pair of distinct edges receive at least three colours. That is, every bichromatic subgraph has at most one edge. The \emph{harmonious chromatic number} \hcn{G} is the minimum number of colours in a harmonious colouring of $G$; see \citep{Edwards-CPC05, Kubale04, HHU-JCMCC03, HHU-JCMCC03a, BBW-JGT89, Edwards97, EM-JGT94, LeeMitchem-JGT87, CE-AJC04,  TV06, Edwards-JGT98, Edwards-JLMS97, KR-JGT94, Edwards-CPC95, Georges-JGT95, MP-DM91, MX-JGT91, HK-SJADM83, Ewards-DAM92, Edwards-DM99, EM-DAM95,  Kubale04, Edwards-CPC96}. 

The following two basic lower bounds on \hcn{G} are well known. Since a vertex and its neighbours all receive distinct colours in a harmonious colouring, 
$$\hcn{G}\geq\Delta+1.$$
Since $G$ has at most $\binom{\hcn{G}}{2}$ edges, 
\begin{equation}
\eqnlabel{HarmEdgesLowerBound}
\hcn{G}>\sqrt{2m}.
\end{equation}

The next bound is new. Observe that a harmonious colouring of $G$ is an oriented colouring for every orientation of $G$. Thus 
\begin{equation}
\eqnlabel{HarmOriented}
\ocn{G}\leq\hcn{G}.
\end{equation}
Many upper bounds on \hcn{G} are known. For example, \citet{MX-JGT91} proved that $\hcn{G}\leq2\Delta\sqrt{n-1}$. Thus \eqnref{HarmOriented} implies the following lemma, which proves the upper bound on \ocn{Q_d} in \eqnref{Hypercube}.

\begin{lemma}
\lemlabel{OrientedUpperBound}
For every graph $G$,
\begin{equation*}
\ocn{G}\;\leq\,2\Delta\sqrt{n-1}\enspace.
\end{equation*} 
\end{lemma}


\section{An Existential Lower Bound}
\seclabel{ExistentialLowerBound}

In this section we construct a graph whose oriented chromatic number is within a constant factor of the upper bound in \lemref{OrientedUpperBound}.

\begin{lemma}
\lemlabel{Construction}
For infinitely many $n$, there is an $n$-vertex $\Delta$-regular graph $G$ with $\ocn{G}=n$ and $\Delta=\sqrt{8n+1}-3$.
\end{lemma}

\begin{proof}
As illustrated in \figref{Construction}, let $G$ be the underlying undirected graph of the oriented graph $D$ with vertex set $\{(i,j):1\leq i<j\leq p\}$ and arcs\\
\hspace*{1em} (a) $(i,j)(i,k)$ whenever $i<j<k$, \\
\hspace*{1em} (b) $(i,j)(k,j)$ whenever $i<k<j$, and\\
\hspace*{1em} (c) $(i,j)(k,i)$ whenever $k<i<j$.

\Figure{Construction}{\includegraphics{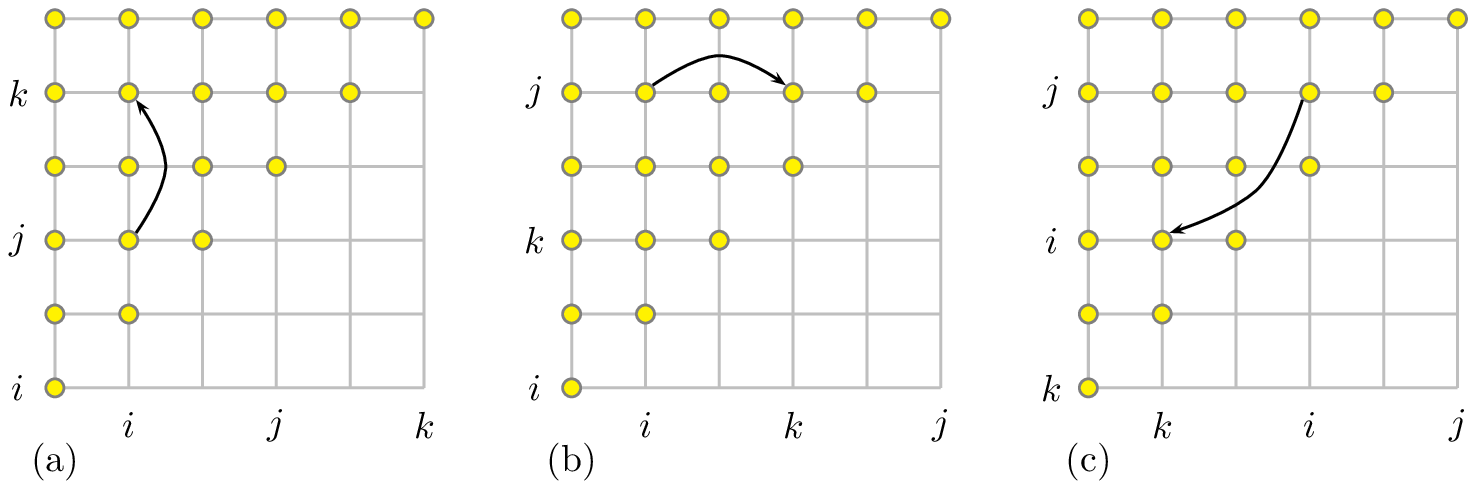}}{Construction of the oriented graph $D$.}

First we compute the degree of each vertex $(i,j)$. Observe that $(i,j)$ has $p-j$ outgoing type-(a) arcs, $j-i-1$ outgoing type-(b) arcs, and $i-1$ outgoing type-(c) arcs. Thus $(i,j)$ has outdegree $(p-j)+(j-i-1)+(i-1)=p-2$. A type-(a) incoming arc at $(i,j)$ is from a vertex $(i,k)$ with $i<k<j$; there are $j-i-1$ such arcs. A type-(b) incoming arc at $(i,j)$ is from a vertex $(k,j)$ with $k<i<j$; there are $i-1$ such arcs. A type-(c) incoming arc at $(i,j)$ is from a vertex $(j,k)$ with $i<j<k$; there are $p-j$ such arcs. Thus $(i,j)$ has outdegree $(j-i+1)+(i-1)+(p-j)=p-2$. Hence $G$ is $\Delta$-regular with $\Delta=2(p-2)=\sqrt{8n+1}-3$.

Suppose on the contrary that $D$ has a directed $2$-cycle $C$. 
If $C$ has a type-(a) arc $(i,j)(i,k)$, then the reverse arc $(i,k)(i,j)$ is also type-(a), implying $j<k$ and $k<j$, which is a contradiction. 
If $C$ has a type-(b) arc $(i,j)(k,j)$, then the reverse arc $(k,j)(i,j)$ is also type-(b), implying $k<j$ and $j<k$, which is a contradiction. 
If $C$ has a type-(c) arc $(i,j)(k,i)$, then the reverse arc is $(k,i)(i,j)$, but there are no arcs of this form. Thus $D$ has no directed $2$-cycle, and indeed $D$ is an oriented graph. 

We claim that $D$ has diameter $2$. Consider two vertices $(i,j)$ and $(k,\ell)$. Then $i<j$ and $k<\ell$, and without loss of generality, $i\leq k$. If $i=k$ and $j<\ell$, then $(i,j)(i,\ell)$ is a type-(a) arc of $D$.
If $i=k$ and $\ell<j$, then $(i,\ell)(i,j)$ is a type-(a) arc of $D$.
Now assume that $i<k$. If $i<k$ and $j=\ell$, then $(i,j)(k,j)$ is a type-(b) arc of $D$. If $i<k$ and $j<\ell$, then $(i,j)(i,\ell)(k,\ell)$ is a type-(ab) path of $D$. Otherwise $i<k$ and $\ell<j$, implying $i<k<\ell<j$, in which case $(i,j)(\ell,j)(k,\ell)$ is a type-(bc) path of $D$. Thus $D$ has diameter $2$, implying $\ocn{D}=n$ by \eqnref{DiameterTwo}.
\end{proof}

It follows from \lemref{Construction} that in any upper bound of the form $\ocn{G}\leq\Oh{\Delta^\alpha n^\beta}$, we must have $\alpha+2\beta\geq2$. In particular with $\beta=\frac{1}{2}$, for the graph $G$ from \lemref{Construction} we have $\ocn{G}=n>\Delta\sqrt{\frac{n}{8}}$. In this sense, the upper bound in \lemref{OrientedUpperBound} is tight up to a constant factor. 

Also note that Moore's bound for the degree/diameter problem implies that $\Delta\geq\sqrt{n-1}$ in every (undirected) graph with diameter $2$. 

\section{A Universal Lower Bound}
\seclabel{UniversalLowerBound}

We now consider universal lower bounds on the oriented chromatic number. \citet{KSZ-JGT97} proved the following\footnote{For completeness we include the proof of \Eqnref{KSZ} by \citet{KSZ-JGT97}. Let $k:=\ocn{G}$. $G$ has less than $k^n$ colourings with $k$ colours, each of which is an oriented colouring of at most $2^{\binom{k}{2}}$ orientations. Thus the number of orientations,  $2^m$, is less than $k^n2^{\binom{k}{2}}$. Thus $m<n\log k+\binom{k}{2}$.}  lower bound for all $G$, which implies \eqnref{DegreeLowerBound}. (Throughout this paper, all logarithms are binary.)\ 
\begin{equation}
\eqnlabel{KSZ}
\binom{\ocn{G}}{2}+n\log\big(\ocn{G}\big)\;\geq\;m\enspace.
\end{equation}
We now reformulate \eqnref{KSZ} for reasonably dense graphs. Say $G$ has average degree $\delta:=\frac{2m}{n}$. Let $t$ be the solution to 
\begin{equation*}
t+\log t\;=\;\delta-\log n\enspace.
\end{equation*}
Note that $0<t<\delta$ and $t\rightarrow\delta$ for $\delta\gg\log n$. (We are not interested in the case $\delta\ll\log n$, when $t$ becomes small.)\ 

\begin{lemma}
\lemlabel{OrientedLowerBound}
For every graph $G$, 
\begin{equation*}
\ocn{G}\;\geq\;\sqrt{nt}\enspace.
\end{equation*}
\end{lemma}

\begin{proof}
Suppose on the contrary that $\ocn{G}<\sqrt{nt}$. By \eqnref{KSZ},
\begin{equation*}
\binom{\sqrt{nt}}{2}+n\log\bracket{\sqrt{nt}}\;>\;m\enspace.
\end{equation*}
Thus $\half nt+\half n\log(nt)>\half \delta n$, implying $t+\log t>\delta-\log n$. This contradiction proves the claim.
\end{proof}

\begin{lemma} 
\lemlabel{OrientedLowerBoundSpecial}
For every graph $G$ with average degree $\delta\geq\log n$, 
\begin{equation*}
\ocn{G}\;\geq\;0.8007...\sqrt{n}\enspace.
\end{equation*}
\end{lemma}

\begin{proof}
\lemref{OrientedLowerBound} implies the claim since $\sqrt{t}\geq0.8007...$
whenever $\delta\geq\log n$.
\end{proof}

For the hypercube, $\delta=\log n$. Thus \lemref{OrientedLowerBoundSpecial} implies the lower bound in \eqnref{Hypercube}. Since \eqnref{KSZ} is proved by a non-constructive counting argument, it would be interesting to construct an orientation $D$ of $Q_d$ with $\ocn{D}\in\Omega(\sqrt{2^d})$; see  
\citep{Matousek-Comb06, FKL-Networks02, EG-IPL89, Baldi88} for results on specific orientations of the hypercube. 

We now refine \lemref{OrientedLowerBound} for graphs that are more dense than hypercubes. 

\begin{lemma} 
\lemlabel{GeneralOrientedLowerBound}
For every graph $G$ with average degree $\delta\geq\log n+(1+\epsilon)\log t$ 
for some $\epsilon>0$,
\begin{equation*}
\ocn{G}\;\geq\;\sqrt{\tfrac{\epsilon}{1+\epsilon}\,(2m-n\log n)}\enspace.
\end{equation*}
\end{lemma}

\noindent(For example, the assumption in \lemref{GeneralOrientedLowerBound} holds if $\delta\geq(2+\epsilon)\log n$.)\ 

\begin{proof}
By the assumption, $t+\log t\geq(1+\epsilon)\log t$ and $t>\epsilon\log t$. Thus \begin{align*}
(1+\epsilon)t\;&>\;\epsilon(t+\log t)\;=\;\epsilon(\delta-\log n)
\hspace*{2em}\text{and}\\
(1+\epsilon)tn\;&>\;\epsilon(\delta n-n\log n)\;=\;\epsilon(2m-n\log n)\enspace.
\end{align*}
Therefore \lemref{OrientedLowerBound} implies that
\begin{equation*}
\ocn{G}
\;\geq\;\sqrt{tn}\;>\;
\sqrt{\tfrac{\epsilon}{1+\epsilon}\,(2m-n\log n)}\enspace.
\end{equation*}
\end{proof}

\lemref{GeneralOrientedLowerBound} says that for sufficiently dense graphs (that is, graphs with super-logarithmic average degree) the lower bound of $\hcn{G}\geq\Omega(\sqrt{m})$ in \eqnref{HarmEdgesLowerBound} also holds for the oriented chromatic number. 

Now suppose that $G$ is $\Delta$-regular for some $\Delta\geq(2+\epsilon)\log n$. Thus \twolemref{OrientedUpperBound}{GeneralOrientedLowerBound} determine \ocn{G} to within a factor of $\Theta(\sqrt{\Delta})$. In particular, 
\begin{equation}
\eqnlabel{RegularBounds}
\sqrt{\tfrac{\epsilon}{2+\epsilon}\,\Delta n}
\;\leq\;
\sqrt{\tfrac{\epsilon}{1+\epsilon}\,(\Delta-\log n)n}
\;\leq\;
\ocn{G}
\;\leq\;
2\Delta\sqrt{n-1}
\enspace.
\end{equation}
The bounds in \eqnref{RegularBounds} are particularly close when $G$ is ($c\log n$)-regular for some constant $c>2$. Then \ocn{G} is between 
$\Omega(\sqrt{n\log n})$ and \Oh{\sqrt{n}\log n}.

\section*{Acknowledgements} 

Thanks to Vida Dujmovi\'c, Andr\'e Raspaud, Bruce Reed, J\'ean-Sebastien Sereni, Ricardo Strausz, and St\'ephan Thomass\'e for stimulating discussions.


\def\soft#1{\leavevmode\setbox0=\hbox{h}\dimen7=\ht0\advance \dimen7
  by-1ex\relax\if t#1\relax\rlap{\raise.6\dimen7
  \hbox{\kern.3ex\char'47}}#1\relax\else\if T#1\relax
  \rlap{\raise.5\dimen7\hbox{\kern1.3ex\char'47}}#1\relax \else\if
  d#1\relax\rlap{\raise.5\dimen7\hbox{\kern.9ex \char'47}}#1\relax\else\if
  D#1\relax\rlap{\raise.5\dimen7 \hbox{\kern1.4ex\char'47}}#1\relax\else\if
  l#1\relax \rlap{\raise.5\dimen7\hbox{\kern.4ex\char'47}}#1\relax \else\if
  L#1\relax\rlap{\raise.5\dimen7\hbox{\kern.7ex
  \char'47}}#1\relax\else\message{accent \string\soft \space #1 not
  defined!}#1\relax\fi\fi\fi\fi\fi\fi} \def\cprime{$'$}
  \def\Dbar{\leavevmode\lower.6ex\hbox to 0pt{\hskip-.23ex \accent"16\hss}D}

\end{document}